\pdfoutput=1
\documentclass{amsart}
\usepackage{amsthm}
\usepackage{amssymb}
\usepackage{amsmath}
\usepackage{graphicx,color}
\usepackage{hyperref}
\theoremstyle{plain}
\newtheorem{theorem}{Theorem}

\newtheorem{fact}{Fact}
\newtheorem{corollary}{Corollary}

\theoremstyle{definition}

\theoremstyle{remark}
\newtheorem{remark}{Remark}

\newcommand{\sign}{\operatorname{sign}}
\newcommand{\id}{\operatorname{id}}
\newcommand{\e}{\operatorname{even}}
\newcommand{\od}{\operatorname{odd}}

\newcommand{\rlaU}{
\,~~\begin{picture}(42,8)
\put(20,8){\circle*{2}}
\put(0,3){\circle{10}}
\put(20,3){\circle{10}}
\put(40,3){\circle{10}}
\put(35.5,3){\vector(-1,0){11}}
\put(4,3){\vector(1,0){11}}
\end{picture}~~
}

\newcommand{\rlatsf}{
\,~~\begin{picture}(42,8)
\put(-2,1){\tiny$3$}
\put(18,1){\tiny$2$}
\put(38,1){\tiny$1$}
\put(20,8){\circle*{2}}
\put(0,3){\circle{10}}
\put(20,3){\circle{10}}
\put(40,3){\circle{10}}
\put(35.5,3){\vector(-1,0){11}}
\put(4,3){\vector(1,0){11}}
\put(28,6){\tiny$+$}
\put(6,6){\tiny $\eta$}
\end{picture}~~
}
\newcommand{\rlaftsMP}{
\,~~\begin{picture}(42,8)
\put(-2,1){\tiny$i$}
\put(18,1){\tiny$k$}
\put(38,1){\tiny$j$}
\put(20,8){\circle*{2}}
\put(0,3){\circle{10}}
\put(20,3){\circle{10}}
\put(40,3){\circle{10}}
\put(35.5,3){\vector(-1,0){11}}
\put(4,3){\vector(1,0){11}}
\put(28,6){\tiny$+$}
\put(6,6){\tiny $-$}
\end{picture}~~
}

\newcommand{\rlatsfINV}{
\,~~\begin{picture}(42,8)
\put(-2,1){\tiny$3$}
\put(18,1){\tiny$2$}
\put(38,1){\tiny$1$}
\put(20,8){\circle*{2}}
\put(0,3){\circle{10}}
\put(20,3){\circle{10}}
\put(40,3){\circle{10}}
\put(35.5,3){\vector(-1,0){11}}
\put(4,3){\vector(1,0){11}}
\put(6,6){\tiny$+$}
\put(28,6){\tiny $\eta$}
\end{picture}~~
}

\newcommand{\rlafst}{
\,~~\begin{picture}(42,8)
\put(-2,1){\tiny$1$}
\put(18,1){\tiny$2$}
\put(38,1){\tiny$3$}
\put(20,8){\circle*{2}}
\put(0,3){\circle{10}}
\put(20,3){\circle{10}}
\put(40,3){\circle{10}}
\put(35.5,3){\vector(-1,0){11}}
\put(4,3){\vector(1,0){11}}
\put(6,6){\tiny$+$}
\put(28,6){\tiny $\eta$}
\end{picture}~~
}
\newcommand{\rlastfPM}{
\,~~\begin{picture}(42,8)
\put(-2,1){\tiny$j$}
\put(18,1){\tiny$k$}
\put(38,1){\tiny$i$}
\put(20,8){\circle*{2}}
\put(0,3){\circle{10}}
\put(20,3){\circle{10}}
\put(40,3){\circle{10}}
\put(35.5,3){\vector(-1,0){11}}
\put(4,3){\vector(1,0){11}}
\put(6,6){\tiny$+$}
\put(28,6){\tiny $-$}
\end{picture}~~
}

\newcommand{\rlaijkU}{
\,~~\begin{picture}(42,8)
\put(-2,1){\tiny$i$}
\put(18,1){\tiny$j$}
\put(38,1){\tiny$k$}
\put(0,8){\circle*{2}}
\put(20,8){\circle*{2}}
\put(40,8){\circle*{2}}
\put(0,3){\circle{10}}
\put(20,3){\circle{10}}
\put(40,3){\circle{10}}
\put(35.5,3){\vector(-1,0){11}}
\put(4,3){\vector(1,0){11}}
\end{picture}~~
}

\newcommand{\rlakjiU}{
\,~~\begin{picture}(42,8)
\put(-2,1){\tiny$k$}
\put(18,1){\tiny$j$}
\put(38,1){\tiny$i$}
\put(0,8){\circle*{2}}
\put(20,8){\circle*{2}}
\put(40,8){\circle*{2}}
\put(0,3){\circle{10}}
\put(20,3){\circle{10}}
\put(40,3){\circle{10}}
\put(35.5,3){\vector(-1,0){11}}
\put(4,3){\vector(1,0){11}}
\end{picture}~~
}

\newcommand{\rlafstINV}{
\,~~\begin{picture}(42,8)
\put(-2,1){\tiny$1$}
\put(18,1){\tiny$2$}
\put(38,1){\tiny$3$}
\put(20,8){\circle*{2}}
\put(0,3){\circle{10}}
\put(20,3){\circle{10}}
\put(40,3){\circle{10}}
\put(35.5,3){\vector(-1,0){11}}
\put(4,3){\vector(1,0){11}}
\put(28,6){\tiny$+$}
\put(6,6){\tiny $\eta$}
\end{picture}~~
}

\newcommand{\rla}{
\,~~\begin{picture}(42,8)
\put(0,8){\circle*{2}}
\put(20,8){\circle*{2}}
\put(40,8){\circle*{2}}
\put(0,3){\circle{10}}
\put(20,3){\circle{10}}
\put(40,3){\circle{10}}
\put(35.5,3){\vector(-1,0){11}}
\put(4,3){\vector(1,0){11}}
\end{picture}~~
}

\newcommand{\rlbU}{
\,~~\begin{picture}(42,8)
\put(20,8){\circle*{2}}
\put(0,3){\circle{10}}
\put(20,3){\circle{10}}
\put(40,3){\circle{10}}
\put(25,3){\vector(1,0){10}}
\put(4,3){\vector(1,0){11}}
\end{picture}~~
}
\newcommand{\rlb}{
\,~~\begin{picture}(42,8)
\put(0,8){\circle*{2}}
\put(20,8){\circle*{2}}
\put(40,8){\circle*{2}}
\put(0,3){\circle{10}}
\put(20,3){\circle{10}}
\put(40,3){\circle{10}}
\put(25,3){\vector(1,0){10}}
\put(4,3){\vector(1,0){11}}
\end{picture}~~
}

\newcommand{\rlbfst}{
\,~~\begin{picture}(42,8)
\put(-2,1){\tiny$1$}
\put(18,1){\tiny$2$}
\put(38,1){\tiny$3$}
\put(6,6){\tiny$+$}
\put(28,6){\tiny $\eta$}
\put(20,8){\circle*{2}}
\put(0,3){\circle{10}}
\put(20,3){\circle{10}}
\put(40,3){\circle{10}}
\put(25,3){\vector(1,0){10}}
\put(4,3){\vector(1,0){11}}
\end{picture}~~
}
\newcommand{\rlbfstPP}{
\,~~\begin{picture}(42,8)
\put(-2,1){\tiny$i$}
\put(18,1){\tiny$j$}
\put(38,1){\tiny$k$}
\put(6,6){\tiny$+$}
\put(28,6){\tiny $+$}
\put(20,8){\circle*{2}}
\put(0,3){\circle{10}}
\put(20,3){\circle{10}}
\put(40,3){\circle{10}}
\put(25,3){\vector(1,0){10}}
\put(4,3){\vector(1,0){11}}
\end{picture}~~
}

\newcommand{\rlbijkU}{
\,~~\begin{picture}(42,8)
\put(-2,1){\tiny$i$}
\put(18,1){\tiny$j$}
\put(38,1){\tiny$k$}
\put(0,8){\circle*{2}}
\put(20,8){\circle*{2}}
\put(40,8){\circle*{2}}
\put(0,3){\circle{10}}
\put(20,3){\circle{10}}
\put(40,3){\circle{10}}
\put(25,3){\vector(1,0){10}}
\put(4,3){\vector(1,0){11}}
\end{picture}~~
}

\newcommand{\rlbfstINV}{
\,~~\begin{picture}(42,8)
\put(-2,1){\tiny$1$}
\put(18,1){\tiny$2$}
\put(38,1){\tiny$3$}
\put(28,6){\tiny$+$}
\put(6,6){\tiny $\eta$}
\put(20,8){\circle*{2}}
\put(0,3){\circle{10}}
\put(20,3){\circle{10}}
\put(40,3){\circle{10}}
\put(25,3){\vector(1,0){10}}
\put(4,3){\vector(1,0){11}}
\end{picture}~~
}

\newcommand{\rlbtsf}{
\,~~\begin{picture}(42,8)
\put(-2,1){\tiny$3$}
\put(18,1){\tiny$2$}
\put(38,1){\tiny$1$}
\put(6,6){\tiny$+$}
\put(28,6){\tiny $\eta$}
\put(20,8){\circle*{2}}
\put(0,3){\circle{10}}
\put(20,3){\circle{10}}
\put(40,3){\circle{10}}
\put(25,3){\vector(1,0){10}}
\put(4,3){\vector(1,0){11}}
\end{picture}~~
}
\newcommand{\rlbtsfINV}{
\,~~\begin{picture}(42,8)
\put(-2,1){\tiny$3$}
\put(18,1){\tiny$2$}
\put(38,1){\tiny$1$}
\put(28,6){\tiny$+$}
\put(6,6){\tiny $\eta$}
\put(20,8){\circle*{2}}
\put(0,3){\circle{10}}
\put(20,3){\circle{10}}
\put(40,3){\circle{10}}
\put(25,3){\vector(1,0){10}}
\put(4,3){\vector(1,0){11}}
\end{picture}~~
}

\newcommand{\rlcU}{
\,~~\begin{picture}(42,8)
\put(20,8){\circle*{2}}
\put(0,3){\circle{10}}
\put(20,3){\circle{10}}
\put(40,3){\circle{10}}
\put(35.5,3){\vector(-1,0){11}}
\put(15,3){\vector(-1,0){11}}
\end{picture}~~
}
\newcommand{\rlcfstINV}{
\,~~\begin{picture}(42,8)
\put(-2,1){\tiny$1$}
\put(18,1){\tiny$2$}
\put(38,1){\tiny$3$}
\put(6,6){\tiny$\eta$}
\put(28,6){\tiny $+$}
\put(20,8){\circle*{2}}
\put(0,3){\circle{10}}
\put(20,3){\circle{10}}
\put(40,3){\circle{10}}
\put(35.5,3){\vector(-1,0){11}}
\put(15,3){\vector(-1,0){11}}
\end{picture}~~
}
\newcommand{\rlcfst}{
\,~~\begin{picture}(42,8)
\put(-2,1){\tiny$1$}
\put(18,1){\tiny$2$}
\put(38,1){\tiny$3$}
\put(28,6){\tiny$\eta$}
\put(6,6){\tiny $+$}
\put(20,8){\circle*{2}}
\put(0,3){\circle{10}}
\put(20,3){\circle{10}}
\put(40,3){\circle{10}}
\put(35.5,3){\vector(-1,0){11}}
\put(15,3){\vector(-1,0){11}}
\end{picture}~~
}
\newcommand{\rlckjiU}{
\,~~\begin{picture}(42,8)
\put(-2,1){\tiny$k$}
\put(18,1){\tiny$j$}
\put(38,1){\tiny$i$}
\put(0,8){\circle*{2}}
\put(20,8){\circle*{2}}
\put(40,8){\circle*{2}}
\put(0,3){\circle{10}}
\put(20,3){\circle{10}}
\put(40,3){\circle{10}}
\put(35.5,3){\vector(-1,0){11}}
\put(15,3){\vector(-1,0){11}}
\end{picture}~~
}
\newcommand{\rlcijkU}{
\,~~\begin{picture}(42,8)
\put(-2,1){\tiny$i$}
\put(18,1){\tiny$j$}
\put(38,1){\tiny$k$}
\put(0,8){\circle*{2}}
\put(20,8){\circle*{2}}
\put(40,8){\circle*{2}}
\put(0,3){\circle{10}}
\put(20,3){\circle{10}}
\put(40,3){\circle{10}}
\put(35.5,3){\vector(-1,0){11}}
\put(15,3){\vector(-1,0){11}}
\end{picture}~~
}

\newcommand{\rlctsfINV}{
\,~~\begin{picture}(42,8)
\put(-2,1){\tiny$3$}
\put(18,1){\tiny$2$}
\put(38,1){\tiny$1$}
\put(6,6){\tiny$\eta$}
\put(28,6){\tiny $+$}
\put(20,8){\circle*{2}}
\put(0,3){\circle{10}}
\put(20,3){\circle{10}}
\put(40,3){\circle{10}}
\put(35.5,3){\vector(-1,0){11}}
\put(15,3){\vector(-1,0){11}}
\end{picture}~~
}

\newcommand{\rlctsfINVPP}{
\,~~\begin{picture}(42,8)
\put(-2,1){\tiny$k$}
\put(18,1){\tiny$j$}
\put(38,1){\tiny$i$}
\put(6,6){\tiny$+$}
\put(28,6){\tiny $+$}
\put(20,8){\circle*{2}}
\put(0,3){\circle{10}}
\put(20,3){\circle{10}}
\put(40,3){\circle{10}}
\put(35.5,3){\vector(-1,0){11}}
\put(15,3){\vector(-1,0){11}}
\end{picture}~~
}

\newcommand{\rlctsf}{
\,~~\begin{picture}(42,8)
\put(-2,1){\tiny$3$}
\put(18,1){\tiny$2$}
\put(38,1){\tiny$1$}
\put(28,6){\tiny$\eta$}
\put(6,6){\tiny $+$}
\put(20,8){\circle*{2}}
\put(0,3){\circle{10}}
\put(20,3){\circle{10}}
\put(40,3){\circle{10}}
\put(35.5,3){\vector(-1,0){11}}
\put(15,3){\vector(-1,0){11}}
\end{picture}~~
}

\newcommand{\rlc}{
\,~~\begin{picture}(42,8)
\put(0,8){\circle*{2}}
\put(20,8){\circle*{2}}
\put(40,8){\circle*{2}}
\put(0,3){\circle{10}}
\put(20,3){\circle{10}}
\put(40,3){\circle{10}}
\put(35.5,3){\vector(-1,0){11}}
\put(15,3){\vector(-1,0){11}}
\end{picture}~~
}

\newcommand{\rldU}{
\,~~\begin{picture}(42,8)
\put(20,8){\circle*{2}}
\put(0,3){\circle{10}}
\put(20,3){\circle{10}}
\put(40,3){\circle{10}}
\put(15,3){\vector(-1,0){11}}
\put(25,3){\vector(1,0){10}}
\end{picture}~~
}
\newcommand{\rldfst}{
\,~~\begin{picture}(42,8)
\put(-2,1){\tiny$1$}
\put(18,1){\tiny$2$}
\put(38,1){\tiny$3$}
\put(0,8){\circle*{2}}
\put(20,8){\circle*{2}}
\put(40,8){\circle*{2}}
\put(0,3){\circle{10}}
\put(20,3){\circle{10}}
\put(40,3){\circle{10}}
\put(15,3){\vector(-1,0){11}}
\put(25,3){\vector(1,0){10}}
\put(6,6){\tiny$+$}
\put(28,6){\tiny $\eta$}
\end{picture}~~
}
\newcommand{\rldkjiU}{
\,~~\begin{picture}(42,8)
\put(-2,1){\tiny$k$}
\put(18,1){\tiny$j$}
\put(38,1){\tiny$i$}
\put(0,8){\circle*{2}}
\put(20,8){\circle*{2}}
\put(40,8){\circle*{2}}
\put(0,3){\circle{10}}
\put(20,3){\circle{10}}
\put(40,3){\circle{10}}
\put(15,3){\vector(-1,0){11}}
\put(25,3){\vector(1,0){10}}
\end{picture}~~
}

\newcommand{\rldijkU}{
\,~~\begin{picture}(42,8)
\put(-2,1){\tiny$i$}
\put(18,1){\tiny$j$}
\put(38,1){\tiny$k$}
\put(0,8){\circle*{2}}
\put(20,8){\circle*{2}}
\put(40,8){\circle*{2}}
\put(0,3){\circle{10}}
\put(20,3){\circle{10}}
\put(40,3){\circle{10}}
\put(15,3){\vector(-1,0){11}}
\put(25,3){\vector(1,0){10}}
\end{picture}~~
}

\newcommand{\rldfstINV}{
\,~~\begin{picture}(42,8)
\put(-2,1){\tiny$1$}
\put(18,1){\tiny$2$}
\put(38,1){\tiny$3$}
\put(0,8){\circle*{2}}
\put(20,8){\circle*{2}}
\put(40,8){\circle*{2}}
\put(0,3){\circle{10}}
\put(20,3){\circle{10}}
\put(40,3){\circle{10}}
\put(15,3){\vector(-1,0){11}}
\put(25,3){\vector(1,0){10}}
\put(28,6){\tiny$+$}
\put(8,6){\tiny $\eta$}
\end{picture}~~
}

\newcommand{\rldtsf}{
\,~~\begin{picture}(42,8)
\put(-2,1){\tiny$3$}
\put(18,1){\tiny$2$}
\put(38,1){\tiny$1$}
\put(0,8){\circle*{2}}
\put(20,8){\circle*{2}}
\put(40,8){\circle*{2}}
\put(0,3){\circle{10}}
\put(20,3){\circle{10}}
\put(40,3){\circle{10}}
\put(15,3){\vector(-1,0){11}}
\put(25,3){\vector(1,0){10}}
\put(6,6){\tiny$+$}
\put(28,6){\tiny $\eta$}
\end{picture}~~
}
\newcommand{\rldtfsMP}{
\,~~\begin{picture}(42,8)
\put(-2,1){\tiny$k$}
\put(18,1){\tiny$i$}
\put(38,1){\tiny$j$}
\put(0,8){\circle*{2}}
\put(20,8){\circle*{2}}
\put(40,8){\circle*{2}}
\put(0,3){\circle{10}}
\put(20,3){\circle{10}}
\put(40,3){\circle{10}}
\put(15,3){\vector(-1,0){11}}
\put(25,3){\vector(1,0){10}}
\put(6,6){\tiny$-$}
\put(28,6){\tiny $+$}
\end{picture}~~
}

\newcommand{\rldtsfINV}{
\,~~\begin{picture}(42,8)
\put(-2,1){\tiny$3$}
\put(18,1){\tiny$2$}
\put(38,1){\tiny$1$}
\put(0,8){\circle*{2}}
\put(20,8){\circle*{2}}
\put(40,8){\circle*{2}}
\put(0,3){\circle{10}}
\put(20,3){\circle{10}}
\put(40,3){\circle{10}}
\put(15,3){\vector(-1,0){11}}
\put(25,3){\vector(1,0){10}}
\put(28,6){\tiny$+$}
\put(8,6){\tiny $\eta$}
\end{picture}~~
}
\newcommand{\rldsftINVPM}{
\,~~\begin{picture}(42,8)
\put(-2,1){\tiny$j$}
\put(18,1){\tiny$i$}
\put(38,1){\tiny$k$}
\put(0,8){\circle*{2}}
\put(20,8){\circle*{2}}
\put(40,8){\circle*{2}}
\put(0,3){\circle{10}}
\put(20,3){\circle{10}}
\put(40,3){\circle{10}}
\put(15,3){\vector(-1,0){11}}
\put(25,3){\vector(1,0){10}}
\put(28,6){\tiny$-$}
\put(6,6){\tiny $+$}
\end{picture}~~
}

\newcommand{\rld}{
\,~~\begin{picture}(42,8)
\put(0,8){\circle*{2}}
\put(20,8){\circle*{2}}
\put(40,8){\circle*{2}}
\put(0,3){\circle{10}}
\put(20,3){\circle{10}}
\put(40,3){\circle{10}}
\put(15,3){\vector(-1,0){11}}
\put(25,3){\vector(1,0){10}}
\end{picture}~~
}

\begin{document}
\title[Milnor's triple linking number and Gauss diagrams]{
Milnor's triple linking number 
and Gauss diagram formulas of 3-bouquet graphs
}
\author{Noboru Ito}
\address{
National Institute of Technology, Ibaraki College, 866 Nakane Hitachinaka, Ibaraki 312-8508, Japan
}
\email{nito@gm.ibaraki-ct.ac.jp}
\author{Natsumi Oyamaguchi}
\address{Department of teacher education, 1-1 Daigaku-cho, Yachiyo City, Shumei University, Chiba 276-0003, Japan}
\email{p-oyamaguchi@mailg.shumei-u.ac.jp}
\date{May 28, 2022}
\maketitle
\begin{abstract} 
In this paper, we introduce two functions such that the subtraction corresponds to the Milnor's triple linking number; the addition obtains a new integer-valued link homotopy invariant of $3$-component links.     
We also have found a series of  integer-valued invariants derived from four terms whose sum equals the Milnor's triple linking number.    We apply this structure to give invariants of $3$-bouquet graphs.  
\end{abstract}

\section{Introduction}
%
%
Bouquet graphs are elementary topological objects which have been well studied.  However, explicit Gauss diagram formulas \cite{Ostlund2004} of $3$-bouquet graphs have not been very few or may be unknown.   

Let us consider flat vertex isotopy classes of $3$-bouquet graphs.   Since any flat vertex isotopy preserves the cyclic order of edges connecting to the flat vertex, we choose a cyclic order and fix it.    
It is graphically explained by the next paragrph.   

\begin{figure}[b]
\includegraphics[width=12cm]{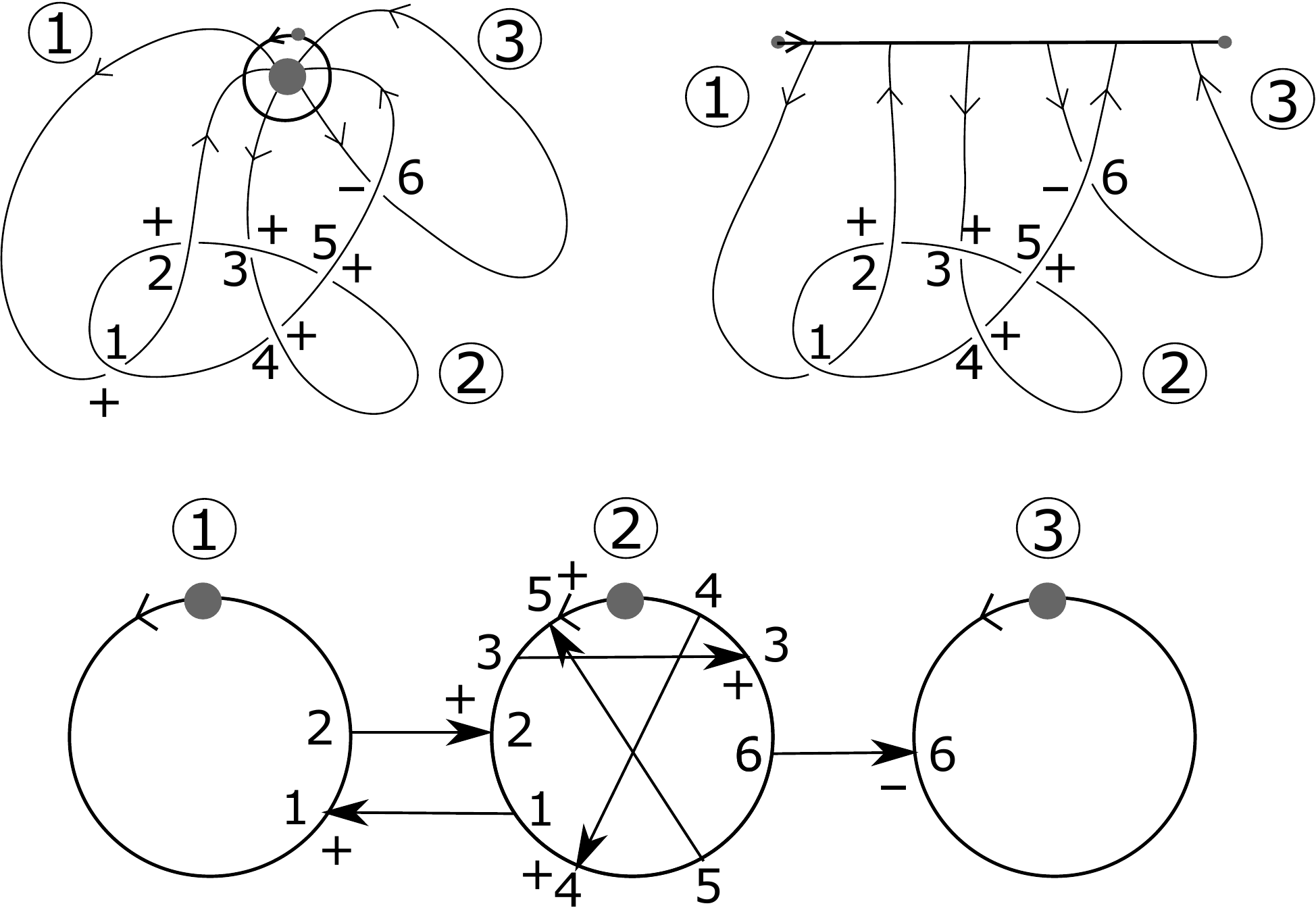}
\caption{A 3-bouquet graph (upper left); a choice of oriented tangles (upper right); its corresponding Gauss diagram (lower line).  
}\label{2strings}
\end{figure}

We firstly take a small disk which center is the flat vertex and assign the fixed cyclic order to  intersections between edges and the boundary of the disk.    It means that we select a base point on the boundary of the disk $d$ (Fig.~\ref{2strings}).   From the base point, we read the endpoints on the boundary by the cyclic order and have the Gauss word: $p_1 p_2 p_3 p_4 p_5 p_6$.  Then for each pair $p_i, p_j$ ($i < j$) belonging to the same component, we orient the component by setting that $p_i$ is the starting point and $p_j$ is the end point.      
Since the boundary $\partial d$ $=$ $S^1$, the isomorphism $S^1 \setminus $ $\{$ a base point $\}$ $\to$ $\mathbb{R}$ induces a mapping from a flat vertex isotopy classes of a $3$-bouquet graph  with a base point on $\partial d$ to a (6, 0)-tangle.      

Throughout this paper, we call the flat vertex isotopy classes of $3$-bouquet graphs equipped with a base point as above \emph{based} flat vertex isotopy classes.  
In the following, every notation obeys \"{O}stlund \cite{Ostlund2004} including Gauss diagram formulas and arrow diagrams.  
 
\begin{theorem}\label{thmUSMilnor}
Let $G_k$ be a Gauss diagram of an ordered three-component link $k$, where $k_1$, $k_2$, $k_3$ are the components of $k$.  Let $\sigma$ $=$ $\left(\begin{matrix} 1&2&3 \\ i & j & k \end{matrix}\right)$; indices $i, j,$ and $k$ are assigned with three circles.    
Let   
\[
P_{\e} (k) = \frac{1}{6} \sum_{\sigma : \e} \langle 2~\rlbijkU + 2~\rlcijkU + \rlaijkU + \rldijkU , G_k \rangle 
\]
and 
\[
P_{\od} (k) = \frac{1}{6} \sum_{\sigma : \od} \langle 2~\rlbijkU + 2~\rlcijkU + \rlaijkU + \rldijkU , G_k \rangle.    
\]
Then  
$
P_{\e} + P_{\od}  
$ 
is an integer-valued base-point-free link homotopy invariant.
In comparison, 
$
P_{\e} - P_{\od}  \mod \gcd(lk(k_2, k_3), lk(k_1, k_3), lk(k_1, k_2)) 
$ 
is the Milnor's triple linking number that is torsion-valued base-point-free.  
\end{theorem}
\begin{remark}\label{SixInteger}
Note that $P_{\e} \pm P_{\od}$ $\in \mathbb{Z}$.  This is because the value of links with no crossings is obviously zero;  the difference of values by applying  a single crossing change is a multiple $6$ \footnote{Dr. Nakagane told NI this way using crossing change.}.   
\end{remark}
\begin{corollary}
Let $t \in \mathbb{Q}$ and $\hat{\mu}$ $=$ $P_{\e} 
+ P_{\od}$.  
Then 
\[
(1-t) \mu + t \hat{\mu}
\]
is link homotopy invariant for links with the fixed base points.  If $t=1$, it is the base-point-free invariant.  
\end{corollary}
\begin{theorem}\label{VersionMilnor}
Let 
\[
P_1 (k) = \sum_{\sigma \in S_3} \sign(\sigma) \langle ~\rlbijkU + \rlcijkU + \rlaijkU, G_k \rangle,
\]
and 
\[
P_2 (k) = \sum_{\sigma \in S_3} \sign(\sigma) \langle ~\rlbijkU + \rlcijkU + \rldijkU, G_k \rangle.   
\] 
Then $P_1$, $P_2$, $P_{\e}+P_{\od}$, and $\mu_{123} = P_{\e}-P_{\od}$ are invariants of based flat vertex isotopy classes of $3$-bouquet graphs and satisfy  
$\mu_{123} (k) = \frac{1}{6} (P_1 (k) + P_2 (k))$.    
We also have that
\[
P_*  \mod \gcd(2lk(k_2, k_3), 2lk(k_1, k_3), 2lk(k_1, k_2))~~~(*=1, 2)
\]
is a base-point-free link homotopy invariant.  
\end{theorem}

\begin{theorem}\label{USthm}
Let $G_k$ be a Gauss diagram of an ordered three-component link $k$.   Let $\sigma$ $=$ $\left(\begin{matrix} 1&2&3 \\ i & j & k \end{matrix}\right)$.   
The $18$ functions
\[
Q^{\sigma}_1 (k) = \langle ~\rlbijkU + \rlckjiU, G_k \rangle, 
\]
\[
Q^{\sigma}_2 (k) = \langle ~\rldijkU + \rldkjiU, G_k \rangle, 
\]
and
\[
Q^{\sigma}_3 (k) = \langle ~\rlaijkU + \rlakjiU, G_k \rangle
\]
are base-point-free link homotopy invariants, each of which is   independent of the Milnor invariant $\mu_{123}$ for $k$.     The $18$ functions are also invariants of based flat vertex isotopy classes of $3$-bouquet graphs. 

If $n \neq n'$, 
\[
(Q^{\id}_n, Q^{(123)}_n, Q^{(132)}_n, Q^{(23)}_n, Q^{(12)}_n, Q^{(13)}_n) \neq  (Q^{\id}_{n'}, Q^{(123)}_{n'}, Q^{(132)}_{n'}, Q^{(23)}_{n'}, Q^{(12)}_{n'}, Q^{(13)}_{n'}).  
\]


\end{theorem}
\begin{figure}
\includegraphics[width=5cm]{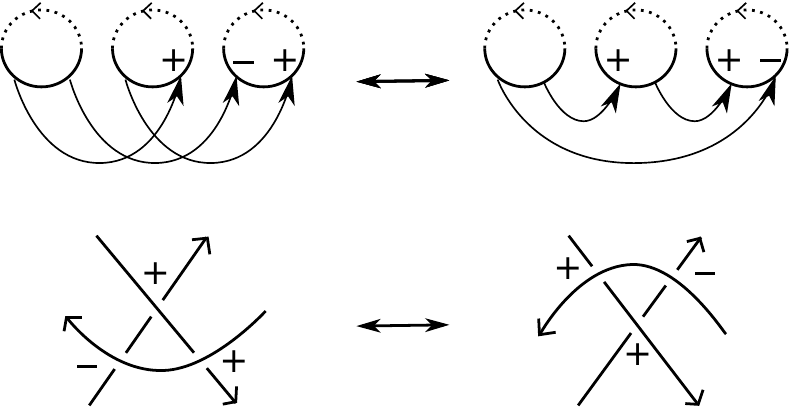}
\caption{A third Reidemeister move with respect to 3-components and the corresponding Gauss diagrams.  We use the move $\Omega_{III+-+3}$   
 included in \cite[Table~1]{Ostlund2004}.}\label{reide&relator}
\end{figure}
\section{Proofs of Theorems~\ref{thmUSMilnor}-- \ref{USthm}}
In this section, every notation of Gauss diagram formulas $\langle \cdot, \cdot  \rangle$ obeys the paper  \cite{Ostlund2004}  (Tables~\ref{baseOver}--\ref{InvIII}).   We use the notion of \emph{diagram fragments} in the  paper \cite[Section~4.4]{Ostlund2004} of \"{O}stlund.   
\subsection{Invariance under changes of base points}\label{Section:base}  
Note that the four equalities hold:
\begin{align*}
&\langle ~~\rlb, \cdot \rangle =   \langle ~~\rlbU, \cdot \rangle, \\
& \langle ~~\rlc, \cdot \rangle  =  \langle ~~\rlcU, \cdot \rangle, \\
&\langle ~~\rld, \cdot \rangle = \langle ~~\rldU, \cdot \rangle,~{\textrm{and}}~\\
&\langle ~~\rla, \cdot \rangle =  \langle ~~\rlaU, \cdot \rangle.
\end{align*}
Then it is sufficient to focus on two kinds of the base point 
moves

\[\begin{matrix}
 \begin{picture}(0,0)
\put(0,11){$\epsilon$}
\qbezier(-10,0)(0,10)(-10,20)
\qbezier(15,0)(5,10)(15,20)
\put(-5,10){\vector(4,-1){15}}
\put(10,10){\circle*{2}}
\end{picture} & \qquad
\begin{picture}(0,0)
\put(-5,8){$\to$}
\end{picture} & \qquad
 \begin{picture}(0,0)
 \put(0,14){$\epsilon$}
\qbezier(-10,0)(0,10)(-10,20)
\qbezier(15,0)(5,10)(15,20)
\put(-5,10){\vector(4,1){15}}
\put(10,10){\circle*{2}}
\end{picture}
\end{matrix} \qquad, \qquad  
\begin{matrix}
\begin{picture}(0,0)
 \put(0,10){$\epsilon$}
\qbezier(-10,0)(0,10)(-10,20)
\qbezier(15,0)(5,10)(15,20)
\put(11,5){\vector(-4,1){16}}
\put(10,10){\circle*{2}}
\end{picture} & \qquad
\begin{picture}(0,0)
\put(-5,8){$\to$}
\end{picture} & \qquad
 \begin{picture}(0,0)
 \put(0,14){$\epsilon$}
\qbezier(-10,0)(0,10)(-10,20)
\qbezier(15,0)(5,10)(15,20)
\put(11,14){\vector(-4,-1){16}}
\put(10,10){\circle*{2}}
\end{picture}
\end{matrix} \qquad
\]
since $-$ sign case is the same up to an overall sign.  
The differences before and after applying base point moves are as in Tables~\ref{baseOver} and \ref{baseUnder}.     
\begin{itemize}
\item For $P_{*}$ ($*=$ $\e$, $\od$) or $Q_*$, after each of these base point moves is applied, the value $\left\langle \cdot, G_k \right\rangle$ differs from the original one by $0$.   
\item For $P_i$ ($i=1, 2$), after each of these base point moves is applied, the value $\left\langle \cdot, G_k \right\rangle$ differs from the original one by $\pm 2 lk(k_1, k_2)$ or $\pm 2 lk(k_2, k_3)$. 
\end{itemize}
Thus, Tables~\ref{baseOver} and \ref{baseUnder} imply the invariances of under the base point moves.  
 
\begin{table}[h!]
\caption{Case~I for permutations $\id$ and $(13)$ (the other cases of permutations are easily recovered seeing them).  Table indicates a  base point move with $+$ sign ($-$ sign case is the same up to an overall sign).}\label{baseOver}
\begin{tabular}{ccccc}\hline \\
\\
Move &&  & & Difference (Right $-$ Left)\\ 
 & \begin{picture}(0,0)
\put(0,11){\tiny$+$}
\qbezier(-10,0)(0,10)(-10,20)
\qbezier(15,0)(5,10)(15,20)
\put(-5,10){\vector(4,-1){15}}
\put(10,10){\circle*{2}}
\end{picture} & 
\begin{picture}(0,0)
\put(-5,8){$\to$}
\end{picture}
 & \begin{picture}(0,0)
\put(0,14){\tiny$+$}
\qbezier(-10,0)(0,10)(-10,20)
\qbezier(15,0)(5,10)(15,20)
\put(-5,10){\vector(4,1){15}}
\put(10,10){\circle*{2}}
\end{picture} & $P_*$-type (upper) \\ &&&& $\hat{\mu}$-, $Q^{\sigma}_*$-type (lower) \\ \hline  &&&& $- 2 lk(k_2, k_3)$  \\
Counted fragment& \rlbfst  & $\to$ & \rlctsfINV &  $0$  \\ \hline &&&& $- 2 lk(k_2, k_3)$  \\
Counted fragment& \rlafst  & $\to$ & \rlatsf &  $0$ \\ \hline  &&&& $2 lk(k_1, k_2)$ \\
Counted fragment& \rlbtsf  & $\to$ & \rlcfstINV & $0$ \\ \hline &&&& $2 lk(k_1, k_2)$ \\
Counted fragment& \rlatsfINV  & $\to$ & \rlafstINV & $0$ \\ \hline \\
\end{tabular}
\end{table} 
\begin{table}[h!]
\caption{Case~I\!I for permutations $\id$ and $(13)$ (the other cases of permutations are easily recovered seeing them).   
Table indicates a  base point move with $+$ sign ($-$ sign case is the same up to an overall sign).}\label{baseUnder}
\begin{tabular}{ccccc}\hline
\\
Move &&  & & Difference (Right $-$ Left) \\
 & \begin{picture}(0,0)
 \put(0,10){$+$}
\qbezier(-10,0)(0,10)(-10,20)
\qbezier(15,0)(5,10)(15,20)
\put(11,5){\vector(-4,1){16}}
\put(10,10){\circle*{2}}
\end{picture} & 
\begin{picture}(0,0)
\put(-5,8){$\to$}
\end{picture}
 & \begin{picture}(0,0)
 \put(0,15){$+$}
\qbezier(-10,0)(0,10)(-10,20)
\qbezier(15,0)(5,10)(15,20)
\put(11,14){\vector(-4,-1){16}}
\put(10,10){\circle*{2}}
\end{picture} 
& $P_*$-type (upper) \\ &&&& $\hat{\mu}$-, $Q^{\sigma}_*$-type (lower) \\ \hline  &&&& $- 2 lk(k_2, k_3)$  \\
Counted fragment & \rlcfst  & $\to$ & \rlbtsfINV & $0$ \\ \hline 
&&&& $- 2 lk(k_2, k_3)$  \\
Counted fragment& \rldfst  & $\to$ & \rldtsfINV & $0$ \\ \hline 
 &&&& $ 2 lk(k_1, k_2)$  \\
Counted fragment & \rlctsf  & $\to$ & \rlbfstINV & $0$ \\ \hline  &&&& $ 2 lk(k_2, k_3)$  \\
Counted fragment& \rldtsf  & $\to$ & \rldfstINV & $0$ \\ \hline \\
\end{tabular}
\end{table} 

\subsection{Invariance under Reidemeister moves}\label{Section:Reide}

\subsubsection{Invariance under Reidemeister moves with respect to one/two component(s)}
We will use the list of Reidemeister moves \cite[Table~1]{Ostlund2004} except for replacing $\Omega_{3+---}$ with $\Omega_{3+-+-}$ as in \cite{Ito2022} \footnote{If you chose the third Reidemeister move $\Omega_{3a}$ of \cite[$\Omega_{3a}$]{Polyak2010}, the corresponding Gauss diagram here is not $\Omega_{3+---}$ but $\Omega_{3+-+-}$.  }, .    
Noting that our formula consisting of four ordered Gauss diagrams $\rlb$, $\rlc$, $\rld$, and $\rla$,   
we have the invariance of Reidemeister moves with respect to one component ($\Omega_{1++}$, $\Omega_{1-+}$, $\Omega_{2+\pm}$, $\Omega_{2-\pm}$, $\Omega_{3+-++}$, $\Omega_{3+---}$) and with respect to two components: $\Omega_{II+\pm}$, $\Omega_{II-\pm}$, $\Omega_{III+-+b}$, $\Omega_{III+-+m}$, $\Omega_{III+-+t}$ \cite[Table~1]{Ostlund2004}.   

\subsubsection{Invariance under $\Omega_{III+-+3}$}
Recall that $\Omega_{III+-+3}$ is as in Fig.~\ref{reide&relator} \cite[Table~1]{Ostlund2004}.  
The complete table of the differences of counted fragments by a single Reidemeister move of type $\Omega_{III+-+3}$ as in Table~\ref{InvIII}.  
The contributions to $P_{\e} \pm P_{\od}$ or $Q^{\sigma}_i$ ($i=1, 2, 3$) do not change before and after applying $\Omega_{III+-+3}$ as in Table~\ref{InvIII}.    We also have that for a given permutation $\sigma$,  the contributions to $P_i (k)$ ($i=1$ or $2$) as in the left column of Table~\ref{InvIII} are canceled out and the right also, respectively.  Thus the invariance of each function holds.  
\begin{table}[h!]
\caption{Invariance under move 
$\Omega_{III+-+3}$.}\label{InvIII}
\begin{tabular}{cccc}\hline
Move&&\includegraphics[width=5cm]{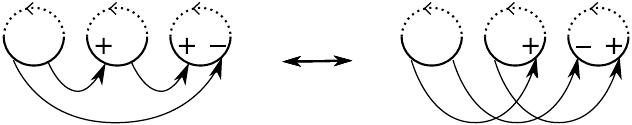}& \\
\\
Counted fragment& \rlbfstPP  & $\longleftrightarrow$ & \rlctsfINVPP \\ \hline \\
Counted fragment& \rlastfPM  & $\longleftrightarrow$ & \rlaftsMP \\ \hline \\
Counted fragment& \rldtfsMP  & $\longleftrightarrow$ & \rldsftINVPM \\ \hline \\
\end{tabular}
\end{table}
\subsection{Identifying our invariants  with Milnor's triple linking number}
Recalling \cite{PolyakViro1994, Ostlund2004} gave  Fact~\ref{PVFormula}, which implies $\mu_{123}$ $=$ $P_{\e} - P_{\od}$.   
\begin{fact}\label{PVFormula}
Let $\sigma$ $=$ $\left(\begin{matrix} 1&2&3 \\ i & j & k \end{matrix}\right)$.   $\mu_{123}(k)$ equals   

\[\frac{1}{6} \sum_{\sigma \in S_3} \sign(\sigma) \left\langle 2~\rlbijkU + 2~\rlcijkU + \rldijkU + \rlaijkU , G_k \right\rangle.  \] 
Then  
$\mu_{123}(k) \mod \gcd( lk (k_2, k_3), lk (k_1, k_3), lk (k_1, k_2))$ is Milnor's triple linking number.    
\end{fact}
\subsection{$\mu_{123}$ ($=P_{\e} - P_{\od}$) and $P_{\e} + P_{\od}$ (or $P_1$, $P_2$) are independent}\label{MainThird}

\begin{figure}[h!]
\includegraphics[width=8cm]{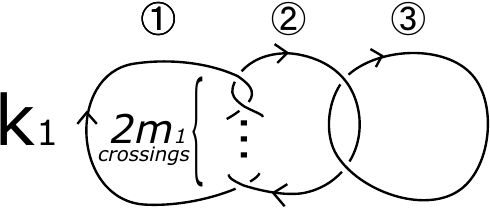}
\caption{Link $k_1$.}\label{202109121}
\end{figure}

\begin{figure}[h!]
\includegraphics[width=8cm]{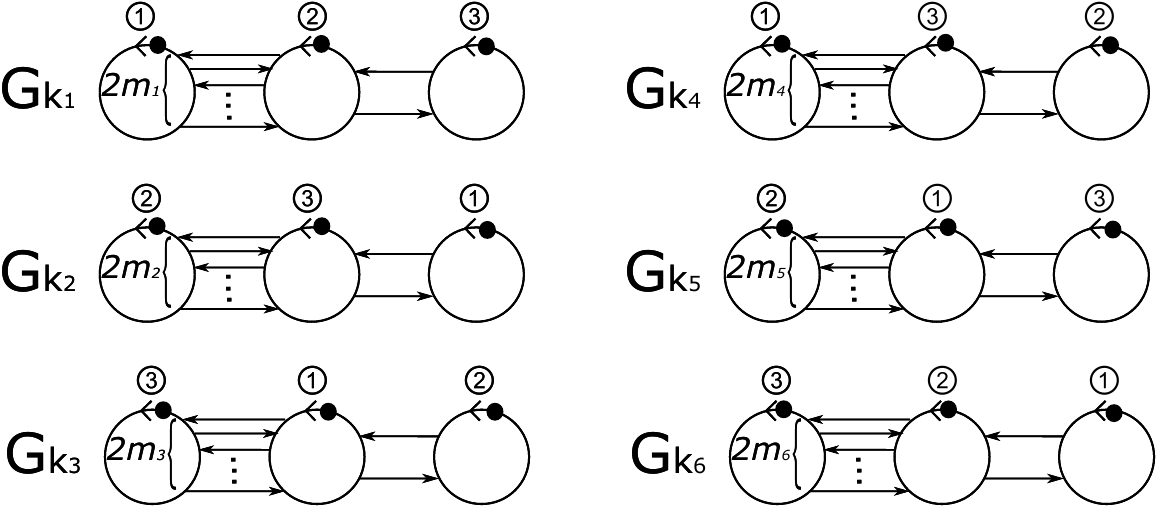}
\caption{Gauss diagrams $G_{k_1}$, $G_{k_2}$,  $G_{k_3}$, $G_{k_4}$, $G_{k_5}$, and   $G_{k_6}$.  Every arrow has the plus sign.}\label{202109122}
\end{figure}

Let $k_1$ be as in Fig.~\ref{202109121} and its Gauss diagram is as in Fig.~\ref{202109122}.  
By definition, $\mu_{123}$ $\equiv$ $0$, whereas $P_{\e}+P_{\od}$ $=$ $m_1$ for a given integer $m_1$.      
By defining similar links $k_i$ ($i=2, 3, 4, 5,$ and $6$), let $G_{k_i}$ be as in Fig.~\ref{202109122}.      
If $m_1$ is odd, Table~\ref{tableP} implies the claim.  
\begin{table}
\caption{Values of $\mu_{123}$ modulo linking numbers and $(P_1, P_2)$ modulo doubled linking numbers for $k_1$ ($m_1$ is odd).}\label{tableP}
\centering
\begin{tabular}{|c|c|} \hline
& $k_1$   \\ \hline
Milnor's $\mu_{123}$ modulo linking numbers & $0$  \\ \hline
$(P_1 (k), P_2 (k))$ & $(1, 1)$   \\ \hline
\end{tabular}
\end{table}
\subsection{Showing integer-valued bouquet graph invariants $P_1$ and $P_2$ are different and they are nontrivial}
Note that a $(6, 0)$-tangle which fixes Gauss diagrams with three base points.   Note also that integer-valued function $P_i$ ($i=1, 2$) is invariant of Reidemeister moves preserving  the base points.  Therefore, we do not need the argument of Section~\ref{Section:base} (i.e. invariance under base point moves are not requested).  
Table~\ref{tableMir} implies the claim.  
\begin{figure}
\includegraphics[width=10cm]{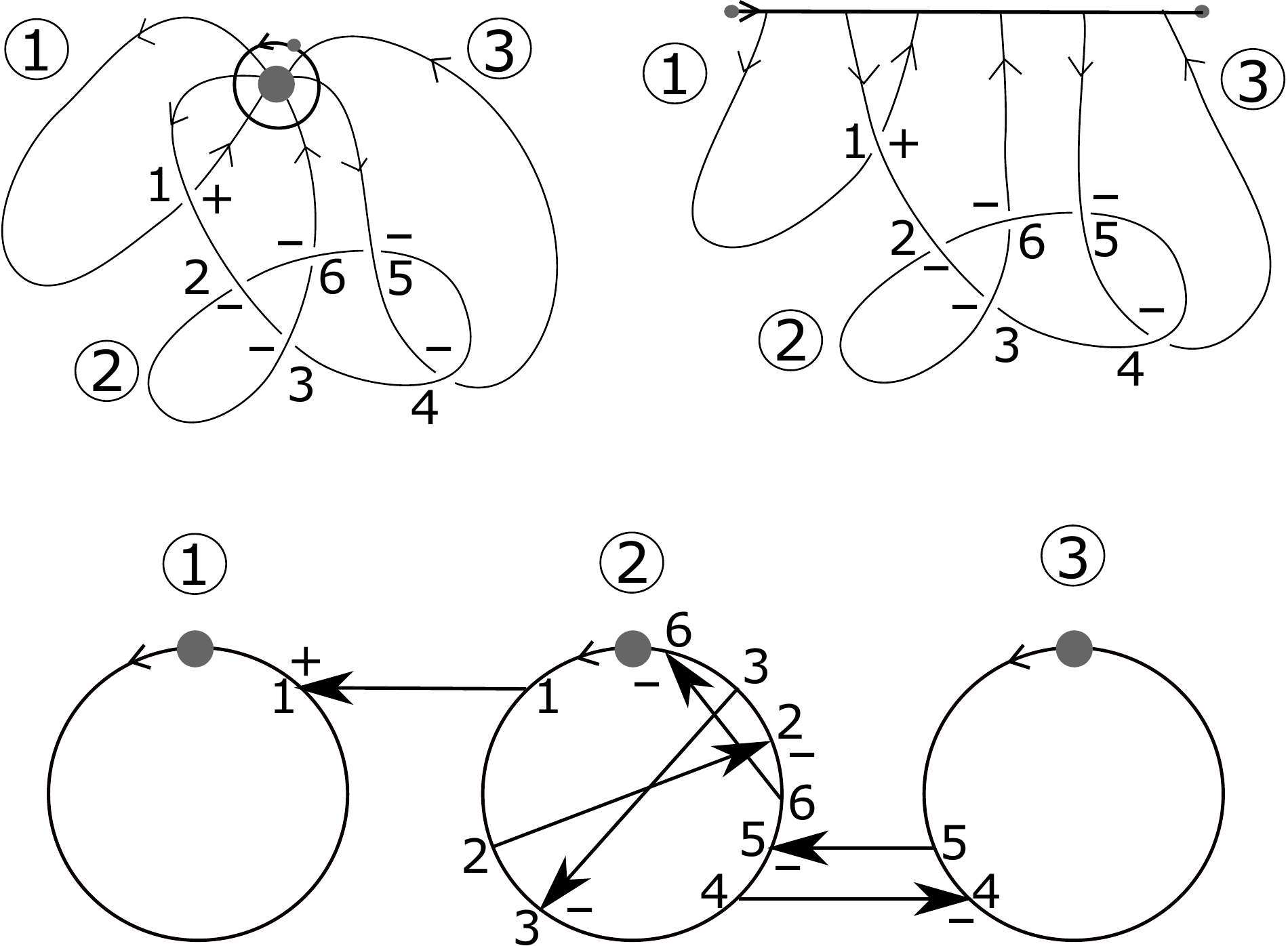}
\caption{The mirror image $b^{\rm{mir}}$.}\label{mirror}
\end{figure}
\begin{table}
\caption{Values of $P_n$ ($n=1, 2$) of the based flat vertex isotopy of a bouquet graph $b$ of Fig.~\ref{2strings} and its mirror image $b^{\rm{mir}}$ (Fig.~\ref{mirror}).}\label{tableMir}
\centering
\begin{tabular}{|c|c|c|} \hline
& $b$ & $b^{\rm{mir}}$   \\ \hline
$\mu_{123} = P_{\e} - P_{\od}$ & $-\frac{1}{2}$ & $-\frac{1}{2}$      \\ \hline
$P_{\e} + P_{\od}$ & $-\frac{1}{2}$ & $-\frac{1}{2}$      \\ \hline
$P_1 (k)$ & $-2$ & $-1$   \\ \hline
$P_2 (k)$ & $-1$ & $-2$   \\ \hline
\end{tabular}
\end{table}


\subsection{Independencies of $Q^{\sigma}_n$}
Let $G_{k_i}$ ($i=1, 2, 3$) be as in Fig.~\ref{202109122}.  

For these links, $\mu_{123}$ $\equiv$ $0$ (Table~\ref{table123}) , which implies the independence of $\mu_{123}$.    
\begin{table}
\caption{Values of $Q^{\sigma}_n$ ($n= 1, 2, 3$) of $k_r$ ($r= 1, 2, 3, 4, 5$, and $6$).}\label{table123}
\centering
\begin{tabular}{|c|c|c|c|c|c|c|} \hline
& $k_1$ & $k_2$ & $k_3$ & $k_4$ & $k_5$ & $k_6$  \\ \hline
Milnor's $\mu_{123}$ modulo linking numbers &0&0&0 &0&0&0  \\ \hline
$Q^{\id}_n (k)$ & $m_1$ & $0$ &  $0$ &0&0&$m_6$  \\ \hline
$Q^{(123)}_n (k)$ & $0$ & $m_2$ & $0$ &$m_4$&0&0  \\ \hline
$Q^{(132)}_n (k)$ & $0$ & $0$ & $m_3$ &0&$m_5$&0  \\ \hline
$Q^{(23)}_n (k)$ & $0$ & $m_2$ & $0$ & $m_4$ &0&0  \\ \hline
$Q^{(12)}_n (k)$ & $0$ & $0$ & $m_3$ & $0$ & $m_5$ & $0$   \\ \hline
$Q^{(13)}_n (k)$ & $m_1$ & $0$ & $0$ & $0$ & $0$ & $m_6$  \\ \hline
\end{tabular}
\end{table}
Table~\ref{tableb} implies the difference $Q^{\id}_n$ and $Q^{\id}_{n'}$ ($n \neq n'$).  
\begin{table}
\caption{Values of $Q^{\operatorname{id}}_1$, ${Q^{\id}}_2$, and ${Q^{\id}}_3$ of the based flat vertex isotopy of a bouquet graph $b$ of Fig.~\ref{2strings} and its mirror image $b^{\rm{mir}}$.}
\label{tableb}
\centering
\begin{tabular}{|c|c|c|} \hline
& $b$ & $b^{\rm{mir}}$ \\ \hline
$Q^{\id}_1$ & $-1$ & $-1$   \\ \hline
$Q^{\id}_2$ & $-1$ & $0$   \\ \hline 
$Q^{\id}_3$ & $0$ & $-1$    \\ \hline
\end{tabular}
\end{table}
\section*{Acknowlegements}
The authors would like to thank Dr.~Keita Nakagane and Dr.~Atsuhiko Mizusawa for their  comments.  The work of NI is partially supported by MEXT KAKENHI Grant Number 20K03604.   

\bibliographystyle{plain}
\bibliography{Ref}
\end{document}